\documentclass[12pt]{article}
\usepackage{amsmath,amsfonts,latexsym,amsthm,amssymb}
\newcommand{\DK}{\mathbb D_{(k)}}
\newcommand{\IK}{\mathbb I_{(k)}}
\newcommand{\K}{\mathcal K}
\newcommand{\Rn}{\mathbb R^n}
\newcommand{\CBF}{\mathcal C\mathcal B\mathcal F}
\newcommand{\wpx}{\widetilde{w}(p,x)}
\newcommand{\Zpx}{\widetilde{Z}(p,x)}
\numberwithin{equation}{section}

\DeclareMathOperator{\I}{Im}

\begin{document}
\newtheorem{teo}{Theorem}
\newtheorem{lem}{Lemma}
\newtheorem{prop}{Proposition}
\newtheorem{cor}{Corollary}
\title{General Fractional Calculus, Evolution Equations,
and Renewal Processes}
\author{Anatoly N. Kochubei
\footnote{Partially supported by the Scientific Program of
the National Academy of Sciences of Ukraine, Project No. 0107U002029,
and by DFG under the project ``Neue Klassen von Evolutionsgleichungen und verwandte Probleme der Spektraltheorie''}
\\ \footnotesize Institute of Mathematics,\\
\footnotesize National Academy of Sciences of Ukraine,\\
\footnotesize Tereshchenkivska 3, Kiev, 01601 Ukraine}
\date{}
\maketitle
\vspace*{2cm}
Running head:\quad  ``General Fractional Calculus''

\vspace*{2cm}
\begin{abstract}
We develop a kind of fractional calculus and theory of relaxation and diffusion equations associated with operators in the time variable, of the form $(\DK u)(t)=\frac{d}{dt}\int\limits_0^tk(t-\tau )u(\tau )\,d\tau -k(t)u(0)$ where $k$ is a nonnegative locally integrable function. Our results are based on the theory of complete Bernstein functions. The solution of the Cauchy problem for the relaxation equation $\DK u=-\lambda u$, $\lambda >0$, proved to be (under some conditions upon $k$) continuous on $[(0,\infty )$ and completely monotone, appears in the description by Meerschaert, Nane, and Vellaisamy of the process $N(E(t))$ as a renewal process. Here $N(t)$ is the Poisson process of intensity $\lambda$, $E(t)$ is an inverse subordinator.
\end{abstract}

\vspace{2cm}


\newpage
\section{INTRODUCTION}

The basic ingredient of the theory of fractional evolution equations \cite{EIK,KST}, the Caputo-Dzhrba\-shyan fractional derivative of order $\alpha \in (0,1)$, has the form
\begin{equation}
(\DK u)(t)=\frac{d}{dt}\int\limits_0^tk(t-\tau )u(\tau )\,d\tau -k(t)u(0)
\end{equation}
where
\begin{equation}
k(t)=\frac{t^{-\alpha }}{\Gamma (1-\alpha )},\quad t>0.
\end{equation}
A recent theory of evolution equations with distributed order derivatives (used for modeling ultraslow relaxation and diffusion processes; see \cite{CGSG,GM,MMGS,MMPG} and references therein) is based on the operator (1.1) with
\begin{equation}
k(t)=\int\limits_0^1\frac{t^{-\alpha }}{\Gamma (1-\alpha )}\, d\rho (\alpha ),\quad t>0,
\end{equation}
where $\rho$ is a Borel measure on $[0,1]$; see \cite{APZ,DB,H,K1,K2,K3,L,MS,UG}.

It is natural to look at a general operator (1.1) and to ask the following question. Under what conditions upon a nonnegative function $k\in L_1^{\text{loc}}(\mathbb R_+)$ does the operator $\DK$ possess a right inverse (a kind of a fractional integral) and produce, as a kind of a fractional derivative, equations of evolution type? The latter means, in particular, that

(A) The Cauchy problem
\begin{equation}
(\DK u)(t)=-\lambda u(t),\quad t>0;\quad u(0)=1,
\end{equation}
where $\lambda >0$, has a unique solution $u_\lambda$, infinitely differentiable for $t>0$ and completely monotone, that is $(-1)^nu_\lambda^{(n)}(t)\ge 0$ for all $t>0$, $n=0,1,2,\ldots$.

(B) The Cauchy problem
\begin{equation}
(\DK w)(t,x)=\Delta w(t,x),\quad t>0,\ x\in \Rn;\quad w(0,x)=w_0(x),
\end{equation}
where $w_0$ is a bounded globally H\"older continuous function, that is $|w_0(\xi )-w_0(\eta )|\le C|\xi -\eta |^\gamma$, $0<\gamma \le 1$, for any $\xi ,\eta \in \Rn$, has a unique bounded solution (the notion of a solution should be defined appropriately). Moreover, the equation in (1.5) possesses a fundamental solution of the Cauchy problem, a kernel with the property of a probability density.

Note that the well-posedness of the Cauchy problem for equations with the operator $\DK$ has been established under much weaker assumptions than those needed for (A) and (B); see \cite{G}.

In the above special cases (A) and (B) are satisfied; see \cite{EIK,K1}. When $k(t)$ has the form (1.2), the function $u_\lambda$ can be expressed via the Mittag-Leffler function $E_\alpha$: $u_\lambda (t)=E_\alpha (-\lambda t^\alpha )$; see \cite{EIK,KST}. In the case (1.3), the asymptotic behavior of $u_\lambda (t)$ is studied in \cite{K1,K3}; for the operator-theoretic meaning of the distributed order derivative and integral see \cite{K2}.

From the point of view of mathematical physics, it is natural to expect the emergence of mathematical theories, in which the relaxation function $u_\lambda (t)$ appears instead of $e^{-\lambda t}$. So far, the first developments of this kind are in the theory of stochastic processes, namely the renewal processes with slowly decaying distribution functions of waiting times; see \cite{MGV,MNV} and references therein.

In this paper we find a class of general operators (1.1) possessing the above evolution generating properties. This class is described in terms of analytic properties of the Laplace transform
\begin{equation}
\K (p)=\int\limits_0^\infty e^{-pt}k(t)\,dt.
\end{equation}
We develop, in particular, a theory of the Cauchy problem (1.5). In contrast to the classical theory of parabolic equations and its analogs known for the cases (1.2) and (1.3), the main technical tools are not the contour integration and explicit estimates, but the theory of complete Bernstein functions \cite{SSV}. Our solution of the Cauchy problem (1.4) leads to an analytic description of general renewal processes constructed (for a slightly more general situation) in \cite{MNV} in terms of the random time change in the classical Poisson process determined by an inverse subordinator.

In Section 2, we give a survey of the results we need about complete Bernstein functions and Stieltjes functions. In Section 3, we introduce and study an analogue, for the general framework, of the fractional integration operator. Section 4 is devoted to the problem (1.4), while Section 5 deals with the problem (1.5).

\section{Complete Bernstein Functions and Stieltjes Functions}

In this section we collect information we need about the classes of functions mentioned in the title. For the detailed exposition see \cite{SSV}.

A real-valued function $f$ on $(0,\infty )$ is called a {\it Bernstein function}, if $f\in C^\infty$, $f(\lambda )\ge 0$ for all $\lambda >0$, and
$$
(-1)^{n-1}f^{(n)}(\lambda )\ge 0\quad \text{for all $n\ge 1,\lambda >0$}.
$$

Equivalently, a function $f:\ (0,\infty )\to \mathbb R$ is a Bernstein function, if and only if
\begin{equation}
f(\lambda )=a+b\lambda +\int\limits_0^\infty \left( 1-e^{-\lambda t}\right)\,\mu (dt)
\end{equation}
where $a,b\ge 0$, and $\mu$ is a Borel measure on $[0,\infty )$, called the {\it L\'evy measure}, such that
\begin{equation}
\int\limits_0^\infty \min (1,t)\,\mu (dt)<\infty .
\end{equation}
The triplet $(a,b,\mu )$ is determined by $f$ uniquely. In particular,
\begin{equation}
a=f(0+),\quad b=\lim\limits_{\lambda \to \infty}\frac{f(\lambda )}{\lambda}.
\end{equation}

A Bernstein function $f$ is said to be a {\it complete Bernstein function}, if its L\'evy measure $\mu$ has a completely monotone density $m(t)$ with respect to the Lebesgue measure, so that (2.1) takes the form
\begin{equation}
f(\lambda )=a+b\lambda +\int\limits_0^\infty \left( 1-e^{-\lambda t}\right) m(t)\,dt
\end{equation}
where, by (2.2),
\begin{equation*}
\int\limits_0^\infty \min (1,t)m(t)\,dt<\infty .
\end{equation*}
Here the complete monotonicity means that $m\in C^\infty (0,\infty )$, $(-1)^nm^{(n)}(t)\ge 0$, $t>0$, for all $n=0,1,2,\ldots$.

Another important class of functions is that of {\it Stieltjes functions}, that is of functions $\varphi$ admitting the integral representation
\begin{equation}
\varphi (\lambda )=\frac{a}{\lambda}+b +\int\limits_0^\infty \frac1{\lambda +t}\sigma (dt)
\end{equation}
where $a,b\ge 0$, $\sigma$ is a Borel measure on $[0,\infty )$, such that
\begin{equation}
\int\limits_0^\infty (1+t)^{-1}\sigma (dt)<\infty .
\end{equation}

Using the identity $(\lambda +t)^{-1}=\int\limits_0^\infty e^{-ts}e^{-\lambda s}ds$ we find from (2.5) that
\begin{equation}
\varphi (\lambda )=\frac{a}{\lambda}+b +\int\limits_0^\infty e^{-\lambda s}g(s)\,ds
\end{equation}
where
\begin{equation}
g(s)=\int\limits_0^\infty e^{-ts}\sigma (dt)
\end{equation}
is a completely monotone function whose Laplace transform exists for any $\lambda >0$.

We will denote the class of complete Bernstein functions by $\CBF$, and the class of Stieltjes functions by $\mathcal S$. The following characterization is proved in \cite{SSV}.

\begin{prop}
Suppose that $f$ is a nonnegative function on $(0,\infty )$. Then the following conditions are equivalent.

\begin{description}
\item[(i)] $f\in \CBF$.

\item[(ii)] The function $\lambda \mapsto \lambda^{-1}f(\lambda )$ is in $\mathcal S$.

\item[(iii)] $f$ has an analytic continuation to the upper half-plane $\mathbb H=\{z\in \mathbb C:\ \I z>0\}$, such that $\I f(z)\ge 0$ for all $z\in \mathbb H$, and there exists the real limit
    \begin{equation}
    f(0+)=\lim\limits_{(0,\infty )\ni \lambda \to 0}f(\lambda ).
\end{equation}

\item[(iv)] $f$ has an analytic continuation to the cut complex plane $\mathbb C\setminus (-\infty ,0]$, such that $\I z\cdot \I f(z)\ge 0$, and there exists the real limit (2.9).

\item[(v)] $f$ has an analytic continuation to $\mathbb H$ given by the expression
    \begin{equation}
    f(z)=a+bz+\int\limits_0^\infty \frac{z}{z+t}\sigma (dt)
    \end{equation}
    where $a,b\ge 0$, and $\sigma$ is a Borel measure on $(0,\infty )$ satisfying (2.6).
\end{description}
\end{prop}

\bigskip
Note that the constants $a,b$ are the same in both the representations (2.4) and (2.10). The density $m(t)$ appearing in the integral representation (2.4) of a function $f\in \CBF$ and the measure $\sigma$ corresponding to the Stieltjes function $\varphi (\lambda )=\lambda^{-1}f(\lambda )$ are connected by the relation
\begin{equation}
m(t)=\int\limits_0^\infty e^{-ts}s\,\sigma (ds).
\end{equation}

The importance of complete Bernstein functions is caused by the following ``nonlinear'' properties \cite{SSV}, quite unusual and having significant applications.

\medskip
\begin{prop}
\begin{description}
\item[(i)] A function $f\not\equiv 0$ is a complete Bernstein function, if and only if $1/f$ is a Stieltjes function.

\item[(ii)] Let $f,f_1,f_2\in \CBF$, $\varphi ,\varphi_1,\varphi_2 \in \mathcal S$. Then $f\circ \varphi \in \mathcal S$, $\varphi \circ f\in \mathcal S$, $f_1\circ f_2\in \CBF$, $\varphi_1\circ \varphi_2 \in \CBF$, $(\lambda +f)^{-1}\in \mathcal S$ for any $\lambda >0$.
\end{description}
\end{prop}

\medskip
It follows from Propositions 1 (ii) and 2 (i) that $0\not\equiv f\in \CBF$, if and only if the function $f^*(\lambda )=\lambda /f(\lambda )$ belongs to $\CBF$. Let us write its representation similar to (2.4),
\begin{equation*}
f^*(\lambda )=a^*+b^*\lambda +\int\limits_0^\infty \left( 1-e^{-\lambda t}\right) m^*(t)\,dt.
\end{equation*}
Then
\begin{equation}
a^*=\begin{cases}
0, & \text{if $a>0$},\\
\frac1{b+\int\limits_0^\infty tm(t)\,dt}, & \text{if $a=0$};
\end{cases}
\end{equation}

\begin{equation}
b^*=\begin{cases}
0, & \text{if $b>0$},\\
\frac1{a+\int\limits_0^\infty m(t)\,dt}, & \text{if $b=0$};
\end{cases}
\end{equation}

\section{Fractional Derivative and Integral}

Throughout the paper (except the uniqueness theorem in Section 5) we make the following assumptions regarding the Laplace transform (1.6) of the function $k$.

\begin{description}
\item[(*)] The Laplace transform (1.6) exists for all $p>0$. The function $\K$ belongs to the Stieltjes class $\mathcal S$, and
    \begin{equation}
    \K (p)\to \infty ,\text{ as $p\to 0$};\quad \K (p)\to 0,\text{ as $p\to \infty$};
    \end{equation}
     \begin{equation}
    p\K (p)\to 0,\text{ as $p\to 0$};\quad p\K (p)\to \infty,\text{ as $p\to \infty$};
    \end{equation}
\end{description}

By Proposition 1, the function $p\mapsto p\K (p)$ is a complete Bernstein function. It follows from (2.3), (3.1), and (3.2) that the integral representations (like (2.5) and (2.7)) of the function $\K$ have the form
$$
\K (p)=\int\limits_0^\infty \frac1{p+t}\sigma (dt)
$$
and
\begin{equation}
\K (p)=\int\limits_0^\infty e^{-ps}g(s)\,ds
\end{equation}
where
$$
g(s)=\int\limits_0^\infty e^{-ts}\sigma (dt),
$$
and the measure $\sigma$ satisfies (2.6). For the function $p\mapsto p\K (p)$ we have
$$
p\K (p)=\int\limits_0^\infty \left( 1-e^{-pt}\right) m(t)\,dt,
$$
and the limit relations from (3.1) and (3.2) show that
\begin{equation}
\int\limits_0^\infty m(t)\,dt=\int\limits_0^\infty tm(t)\,dt=\infty .
\end{equation}

It follows from the uniqueness theorem for the Laplace transform that $g(s)=k(s)$, so that the assumptions $(*)$ imply the representation
\begin{equation}
k(s)=\int\limits_0^\infty e^{-ts}\sigma (dt),\quad 0<s<\infty ,
\end{equation}
so that $k$ is completely monotone.

For each fixed $s\ge 1$, the function $t\mapsto (1+t)e^{-ts}$ is monotone decreasing on $[0,\infty )$, and its value at the origin is 1. It follows from (3.5), (2.6), and the dominated convergence theorem that $k(s)\to 0$, $s\to \infty$.

Note that the conditions $(*)$ are obviously satisfied for the kernel (1.2); under some conditions upon the weight measure $\rho$, they are satisfied also for the case (1.3). On the other hand, given a function $\K$ satisfying $(*)$, one can restore $k$ by the formula (3.5). As a simple example, consider the complete Bernstein function $p\mapsto \log (1+p^\beta )$, $0<\beta <1$ (see Example 15.4.59 in \cite{SSV}), and set $\K (p)=p^{-1}\log (1+p^\beta )$. Then $\K (p)\sim p^{\beta -1}$, as $p\to 0$, and $\K (p)\sim \beta p^{-1}\log p$, as $p\to \infty$, so that the conditions (3.1) and (3.2) are satisfied. The above asymptotic properties are different from those corresponding to the cases (1.2) and (1.3).

By Proposition 2, the function $p\mapsto \dfrac1{p\K (p)}$ belongs to the Stieltjes class. Using (2.12), (2.13), and (3.3) we find its representation similar to (3.3), that is
$$
\dfrac1{p\K (p)} =\int\limits_0^\infty e^{-ps}\varkappa (s)\,ds
$$
where $\varkappa (s)$ is a completely monotone function, $\varkappa (s)\to 0$, as $s\to 0$. Just as in (3.5), we get the representation
\begin{equation}
\varkappa (t)=\int\limits_0^\infty e^{-\lambda t}\eta (d\lambda ),\quad 0<t<\infty ,
\end{equation}
where $\int\limits_0^\infty \dfrac{\eta (d\lambda )}{1+\lambda }<\infty $.

Let us consider the convolution
$$
(k*\varkappa )(t)=\int\limits_0^tk(t-\tau )\varkappa (\tau )\,d\tau .
$$
By the construction of $\varkappa$, the Laplace transform of $k*\varkappa$ equals $\dfrac1p$, so that
\begin{equation}
(k*\varkappa )(t)\equiv 1.
\end{equation}

In other words, $k$ and $\varkappa$ form a pair of Sonine kernels. Such kernels were searched for from the 19th century; see \cite{SC} for a survey. The connection between the complete Bernstein functions and Sonine kernels seems never noticed, but can be useful due to the availability of tables of complete Bernstein functions \cite{SSV}.

Let us study, under the assumptions $(*)$, the generalized fractional differentiation operator $\DK$ of the form (1.1) and the generalized fractional integration operator
\begin{equation}
(\IK f)(t)=\int\limits_0^t\varkappa (t-s)f(s)\,ds.
\end{equation}
The operator $\DK u$ is defined on continuous functions $u$, such that $k*u$ is almost everywhere differentiable, for example, on absolutely continuous functions $u$. The operator $\IK$ is defined on $L_1^{\text{loc}}(0,\infty )$.

The following result extends a well-known property of the Caputo-Dzhrbashyan fractional derivative (see Lemma 2.21 and Lemma 2.22 in \cite{KST}).

\medskip
\begin{teo}
\begin{description}
\item[(i)] If $f$ is a locally bounded measurable function on $(0,\infty )$, then $\DK \IK f=f$.
\item[(ii)] If a function $u$ is absolutely continuous on $[0,\infty )$, then $(\IK \DK u)(t)=u(t)-u(0)$.
\end{description}
\end{teo}

\medskip
{\it Proof}. (i) Let $v=\IK f$. By (3.6) and (3.8),
$$
v(t)=\int\limits_0^t f(s)\,ds\int\limits_0^\infty e^{-\lambda (t-s)}\eta (d\lambda ),
$$
which implies the inequality
\begin{equation}
|v(t)|\le C\int\limits_0^\infty \frac1\lambda \left( 1-e^{-\lambda t}\right)\eta (d\lambda ),\quad 0<t\le T<\infty ,
\end{equation}
(here and below the letter $C$ denotes various positive constants). By (3.9),
$$
|v(t)|\le Ct\int\limits_0^1\eta (d\lambda )+C\int\limits_1^\infty  \left( 1-e^{-\lambda t}\right)\frac{\eta (d\lambda )}\lambda ,
$$
so that $v(t)\to 0$, as $t\to +0$.

By (3.7), we have $k*v=k*\varkappa *f=1*f$, so that $(k*v)(t)=\int\limits_0^tf(\tau )\,d\tau $ is absolutely continuous, and $\DK v=f$.

(ii) We find, using (3.7) twice, that
\begin{multline*}
(\IK \DK u)(t)=\int\limits_0^t\varkappa (t-s)\left[ \frac{d}{ds}\int\limits_0^sk(s-\tau )u(\tau )\,d\tau \right] \,ds -u(0)\\
=-\int\limits_0^t\varkappa (\theta )\left[ \frac{d}{d\theta }\int\limits_0^{t-\theta }k(t-\theta -\tau )u(\tau )\,d\tau \right]
\,d\theta -u(0)=\frac{d}{dt}\int\limits_0^t\varkappa (\theta )
\int\limits_0^{t-\theta }k(t-\theta -\tau )u(\tau )\,d\tau -u(0)\\
=\frac{d}{dt}(\varkappa *k*u)(t)-u(0)=\frac{d}{dt}\int\limits_0^t u(\tau )\,d\tau =u(t)-u(0).\qquad \blacksquare
\end{multline*}

\section{Relaxation Equation and Renewal Processes}

Let us consider the Cauchy problem (1.4).

\medskip
\begin{teo}
Under the assumption $(*)$, the problem (1.4) has a unique solution $u_\lambda (t)$, continuous on $[0,\infty )$, infinitely differentiable and completely monotone on $(0,\infty )$.
\end{teo}

\medskip
{\it Proof}. Applying formally the Laplace transform, we find for the image $\widetilde{u_\lambda }$ of a solution the expression
\begin{equation}
\widetilde{u_\lambda }(p)=\frac{\K (p)}{p\K (p)+\lambda },
\end{equation}
that is $p\widetilde{u_\lambda }(p)=\Phi \circ (p\K (p))$ where $\Phi (z)=\dfrac{z}{z+\lambda }$. By Proposition 1(v), $\Phi$ is a complete Bernstein function. It follows from the assumption $(*)$ and Proposition 2(ii) that $p\mapsto p\widetilde{u_\lambda }(p)$ is a complete Bernstein function. Therefore $\widetilde{u_\lambda }$ is a Stieltjes function, so that by (2.7),
$$
\frac{\K (p)}{p\K (p)+\lambda }=\frac{a}p+b+\int\limits_0^\infty e^{-ps}g(s)\,ds,
$$
where $g$ is completely monotone, possesses the representation (2.8) with the measure satisfying (2.6), and by (2.3), (3.1), and (3.2),
$$
a=\lim\limits_{p\to +0}\frac{p\K (p)}{p\K (p)+\lambda }=0,\quad b=\lim\limits_{p\to \infty }\frac{\K (p)}{p\K (p)+\lambda }=0.
$$

Thus, we have found that
$$
\frac{\K (p)}{p\K (p)+\lambda }=\int\limits_0^\infty e^{-ps}g(s)\,ds,
$$
and the identity (4.1) will hold if we set $u_\lambda =g$. In addition, it follows from (3.2) that
$$
\frac{\K (p)}{p\K (p)+\lambda }\sim \frac1p,\quad p\to \infty .
$$
Taking into account the monotonicity of $g$, we may apply the Karamata-Feller Tauberian theorem \cite{Fe} and find that $u_\lambda (t)\to 1$, as $t\to +0$.

Now it follows from the uniqueness theorem for the Laplace transform that $u_\lambda$ is a solution, in the strong sense, of the problem (1.4). $\qquad \blacksquare$

\medskip
{\it Remark}. It is instructive to see what happens if $\K (p)$ satisfies weaker conditions, for example, if $p\mapsto p\K (p)$ is a Bernstein function, but not a complete Bernstein function. It is seen from comparing (2.1) and (2.4) that the simplest example of this kind is obtained if we take $\K_0(p)=p^{-1}(1-e^{-p})$. Then
$$
\K_0(p)=\int\limits_0^\infty k_0(t)e^{-pt}\,dt,\quad k_0(t)=
\begin{cases}
1, & \text{ if $0\le t\le 1$,}\\
0, & \text{ if $t>1$.}
\end{cases}
$$
It is easy to check that
$$
\left( \mathbb D_{(k_0)}u\right) (t)=\begin{cases}
u(t)-u(t-1), & \text{ if $t>1$;}\\
u(t)-u(0), & \text{ if $0<t\le 1$},
\end{cases}
$$
and a general solution of corresponding problem (1.4) has the form
$$
u(t)=\begin{cases}
1, & \text{ if $t=0$;}\\
(1+\lambda )^{-1}, & \text{ if $0<t\le 1$;}\\
c(1+\lambda )^{-t}, & \text{ if $t>1$},
\end{cases}
$$
$c=\text{const}$. This function is not continuous at the origin. It is continuous at $t=1$, if $c=1$, but even in this case it is not differentiable at $t=1$.

\bigskip
A probabilistic interpretation of Theorem 2 can be given on the basis of the results by Meerschaert, Nane, and Vellaisamy \cite{MNV}. Let $D(t)$ be a subordinator (see \cite{B} for the necessary notions and results from the theory of L\'evy processes),
$$
\mathbf E\left[ e^{-sD(t)}\right] =e^{-t\Psi (s)}
$$
where the Laplace exponent $\Psi$ is a Bernstein function having the representation
$$
\Psi (s)=bs+\int\limits_0^\infty \left( 1-e^{-s\tau }\right) \Phi (d\tau )
$$
with the drift coefficient $b\ge 0$ and the L\'evy measure $\Phi$, such that either $b>0$, or $\Phi ((0,\infty ))=\infty$, or both.

The process $D$ is strictly increasing, thus it possesses an inverse function
$$
E(t)=\inf \{ r>0:\ D(r)>t\} .
$$

Let $N(t)$ be the Poisson process with intensity $\lambda$. It is shown in \cite{MNV} that $N(E(t))$ is a renewal process with independent identically distributed waiting times $J_n$, such that
\begin{equation}
\mathbf P[J_n>t]=\mathbf E\left[ e^{-\lambda E(t)}\right] .
\end{equation}
The Laplace transform of the expression in (4.2) is as follows:
\begin{equation}
\int\limits_0^\infty e^{-st}\mathbf E\left[ e^{-\lambda E(t)}\right]\,dt=\frac{\Psi (s)}{s(\lambda +\Psi (s))}.
\end{equation}

It is interesting to know analytic properties of the function (4.2), and we can find them for the case where $\Psi (s)=s\K (s)$ is a complete Bernstein function (so that $\K$ is a Stieltjes function) satisfying (3.1) and (3.2). In this case the right-hand side of (4.3)
is exactly the Laplace transform of the function $u_\lambda (t)$ described in Theorem 2. Therefore, under these assumptions, the function (4.2) is continuous on $[0,\infty )$ and completely monotone. Earlier such properties were known \cite{MNV} for the cases of stable subordinators and their distributed order extensions. In the latter case they were obtained from a rather complicated explicit integral representation \cite{K1} for $u_\lambda (t)$, so that here we have given a simpler proof based on the theory of complete Bernstein functions.

\section{$\DK$-Heat Equation}

Let us consider the Cauchy problem (1.5). Applying formally the Laplace transform in $t$ to both sides of (1.5) we obtain the following equation for the Laplace transform $\wpx$ of a solution of (1.5):
\begin{equation}
p\K (p)\wpx -\K (p)w_0(x)=\Delta \wpx ,\quad p>0,\ x\in \Rn .
\end{equation}

A bounded function $w(t,x)$ will be called a {\it LT-solution} of the problem (1.5), if $w$ is continuous in $t$ on $[0,\infty )$ uniformly with respect to $x\in \Rn$, $w(0,x)=w_0(x)$, while its Laplace transform $\wpx$ is twice continuously differentiable in $x$, for each $p>0$, and satisfies the equation (5.1).

\medskip
\begin{teo}
Suppose that the assumption $(*)$ is satisfied. There exists such a nonnegative function $Z(t,x)$, $t>0$, $x\in \Rn$, $x\ne 0$, locally integrable in $t$ and infinitely differentiable in $x\ne 0$, that
\begin{equation}
\int\limits_{\Rn }Z(t,x)\,dx=1,\quad t>0,
\end{equation}
and for any bounded globally H\"older continuous function $w_0$, the function
\begin{equation}
w(t,x)=\int\limits_{\Rn }Z(t,x-\xi )w_0(\xi )\,d\xi
\end{equation}
is a LT-solution of the Cauchy problem (1.5).
\end{teo}

\medskip
{\it Proof}. Consider the function
\begin{equation}
g(s,p)=\K (p)e^{-sp\K (p)},\quad s>0,\ p>0.
\end{equation}
Since $p\mapsto p\K (p)$ is a Bernstein function, the function $p\mapsto e^{-sp\K (p)}$ is completely monotone (see conditions for the complete monotonicity in Chapter 13 of \cite{Fe}). By Bernstein's theorem, for each $s\ge 0$, there exists such a probability measure $\mu_s(d\tau )$ that
\begin{equation}
e^{-sp\K (p)}=\int\limits_0^\infty e^{-p\tau }\mu_s(d\tau ).
\end{equation}
The family of measures $\{ \mu_s\}$ is weakly continuous in $s$.

Set
$$
G(s,t)=\int\limits_0^tk(t-\tau )\mu_s(d\tau ).
$$
By (3.3) and (5.5), the Laplace transform in $t$ of the function $G$ coicides with the function $g(s,p)$:
$$
g(s,p)=\int\limits_0^\infty e^{-pt}G(s,t)\,dt.
$$

On the other hand, it is seen from (5.4) that
$$
\int\limits_0^\infty g(s,p)\,ds=\frac1p ,
$$
so that
$$
\int\limits_0^\infty e^{-pt}\,dt \int\limits_0^\infty G(s,t)\,ds=\frac1p ,
$$
which implies the equality
\begin{equation}
\int\limits_0^\infty G(s,t)\,ds=1.
\end{equation}

We define $Z$ by the subordination equality
\begin{equation}
Z(t,x)=\int\limits_0^\infty (4\pi s)^{-n/2}e^{-\frac{|x|^2}{4s}} G(s,t)\,ds,\quad x\ne 0.
\end{equation}
The equality (5.2) follows from (5.6) and properties of the fundamental solution of the classical heat equation. It follows from (5.7) and (5.6) that $Z(t,x)$ is infinitely differentiable in $x\ne 0$.

By the definition of $G$, the Laplace transform in $t$ of the function $Z$,
$$
\Zpx =\int\limits_0^\infty (4\pi s)^{-n/2}e^{-\frac{|x|^2}{4s}}g(s,p)\,ds=\K (p)\int\limits_0^\infty (4\pi s)^{-n/2}e^{-\frac{|x|^2}{4s}}e^{-sp\K (p)}\,ds
$$
exists for $p>0$; by Fubini's theorem, this implies the local integrability of $Z(t,x)$ in $t$. Moreover, using the identity 2.3.16.1 from \cite{PBM} we find that
\begin{equation}
\widetilde{Z}(p,x)=(2\pi )^{-\frac{n}2}|x|^{1-\frac{n}2}\K (p)
(p\K (p))^{\frac{1}2(\frac{n}2-1)}K_{\frac{n}2-1}(|x|\sqrt{p\K
(p)}),\quad x\ne 0,
\end{equation}
where $K_\nu$ is the McDonald function. Note that, by virtue of Proposition 2, $\K (p)\ne 0$ for any $p>0$.

Let us consider the function (5.3) starting from its behavior as $t\to 0$. Using (5.2) we write
$$
w(t,x)-w_0(x)=\int\limits_{\Rn}Z(t,x-\xi )[w_0(\xi )-w_0(x)]\,d\xi ,
$$
so that
$$
|w(t,x)-w_0(x)|\le C\int\limits_0^\infty s^{-n/2}G(s,t)\,ds\int\limits_{\Rn}e^{-\frac{|z|^2}{4s}}|z|^\gamma dz
=C_1\int\limits_0^\infty s^{\gamma /2}G(s,t)\,ds,
$$
and in order to prove that $w(t,x)\to w_0(x)$ (uniformly in $x$), as $t\to 0$, it suffices to show that the function
$$
a(t)=\int\limits_0^\infty s^{\gamma /2}G(s,t)\,ds
$$
tends to 0, as $t\to 0$. Consider its Laplace transform
$$
\widetilde{a}(p)=\int\limits_0^\infty s^{\gamma /2}g(s,p)\,ds
=\K (p)\int\limits_0^\infty s^{\gamma /2}e^{-sp\K (p)}\,ds
=hp^{-1}[p\K (p)]^{-\gamma /2}
$$
where $h>0$.

Since the function $\lambda \mapsto \lambda^{\gamma /2}$ is a complete Bernstein function, we find from Proposition 2 that $p\mapsto [p\K (p)]^{\gamma /2}$ is a complete Bernstein function, the function
$p\mapsto [p\K (p)]^{-\gamma /2}$ is a Stieltjes function possessing a representation of the form (2.7). Thus,
\begin{equation}
[p\K (p)]^{-\gamma /2}=\frac{c}p+d+\int\limits_0^\infty e^{-ps}g_\gamma (s)\,ds
\end{equation}
where $g_\gamma$ is a completely monotone function of the form (2.8) with the measure $\sigma$ satisfying (2.6), $c,d\ge 0$.

To find the constants $c$ and $d$ in (5.9), we can consider the complete Bernstein function
$$
p[p\K (p)]^{-\gamma /2}=c+dp+p\int\limits_0^\infty e^{-ps}g_\gamma (s)\,ds
$$
and use the formulas (2.3). It follows from (3.1) and (3.2) that
$$
c=\lim\limits_{p\to +0}p[p\K (p)]^{-\gamma /2}=\lim\limits_{p\to +0}p^{1-\frac\gamma 2}[\K (p)]^{-\gamma /2}=0,
$$
$$
d=\lim\limits_{p\to \infty }[p\K (p)]^{-\gamma /2}=0.
$$
We have found that
$$
\widetilde{a}(p)=hp^{-1}\int\limits_0^\infty e^{-ps}g_\gamma (s)\,ds
=h\int\limits_0^\infty e^{-pt}\,dt \int\limits_0^tg_\gamma (\tau )\,d\tau ,
$$
so that
$$
a(t)=h\int\limits_0^tg_\gamma (\tau )\,d\tau \longrightarrow 0, \quad \text{as $t\to 0$}.
$$

In order to perform the Laplace transform of the both sides of (5.3), we need estimates of the function (5.8). It is known \cite{BE} that the function $K_\nu$ and all its derivatives decay exponentially at infinity. If $\nu >0$ is not an integer, then
$$
K_\nu (z)=\frac{\pi }{2\sin \nu \pi }\left[ I_{-\nu }(z)-I_\nu (z)\right]
$$
where $I_\nu (z)=z^\nu \varphi_\nu (z)$, and $\varphi_\nu$ is an entire function. Therefore near the origin
\begin{equation}
\left| K_\nu^{(l)}(z)\right| \le C|z|^{-\nu -l},\quad l=0,1,2,\ldots .
\end{equation}

If $\nu$ is a nonnegative integer, then
$$
K_\nu (z)=(-1)^{\nu +1}I_\nu (z)\log \frac{z}2+z^{-\nu }P_\nu (z)+Q_\nu (z)
$$
where $P_\nu$ is a polynomial, $Q_\nu$ is an analytic function on a neighbourhood of the origin, so that the inequality (5.10) remains valid for all natural numbers $\nu$, and
\begin{equation}
\left| K_0^{(l)}(z)\right| \le \begin{cases}
C|\log z|, & \text{ if $l=0$},\\
C|z|^{-l}, & \text{ if $l\ge 1$}.
\end{cases}
\end{equation}

As a result, it follows from (5.8), (5.10), and (5.11) that for each fixed $p$, $j=1,\ldots ,n$,
\begin{equation}
\left| \frac{\partial^l}{\partial x_j^l}\Zpx \right| \le C|x|^{-n+2-l}e^{-a|x|},\quad l=0,1,2,\ldots ,
\end{equation}
$a>0$, if $n\ge 3$. If $n=2$, then (5.12) remains valid for $l\ge 1$, while
\begin{equation}
\left| \Zpx \right| \le C|\log |x||e^{-a|x|},\quad n=2.
\end{equation}
If $n=1$, there is no singularity at $x=0$:
\begin{equation}
\left| \frac{\partial^l}{\partial x^l}\Zpx \right| \le Ce^{-a|x|},\quad n=1.
\end{equation}

Now the Laplace transform in (5.3) is legitimate, so that
$$
\wpx =\int\limits_{\Rn} \widetilde{Z}(p,x-\xi )w_0(\xi )\,d\xi ,\quad p>0,\ x\in \Rn .
$$
The estimates (5.12)-(5.14) justify a single differentiation in spatial variables:
$$
\frac{\partial \wpx }{\partial x_j}=\int\limits_{\Rn} \frac{\partial }{\partial x_j}\widetilde{Z}(p,x-\xi )w_0(\xi )\,d\xi ,\quad j=1,\ldots ,n.
$$
If $n=1$, also the direct second differentiation is possible, thus we will consider the case where $n\ge 2$.

It follows from (5.2) that
\begin{equation}
\int\limits_{\Rn}\Zpx \,dx=\frac1p \quad \text{for all $x$},
\end{equation}
so that
$$
\int\limits_{\Rn}\frac{\partial }{\partial x_j}\Zpx \,dx=0
$$
and
\begin{equation}
\frac{\partial \wpx }{\partial x_j}=\int\limits_{\Rn} \frac{\partial }{\partial x_j}\widetilde{Z}(p,x-\xi )\left[ [w_0(\xi )-w_0(x^0)\right] \,d\xi
\end{equation}
where $x^0\in \Rn$ is an arbitrary fixed point.

Let us divide the domain of integration in (5.16) into two subdomains
$$
\Omega_1=\{ \xi \in \Rn :\ |\xi -x^0|\ge \delta \},\quad \Omega_2=\Rn \setminus \Omega_1,
$$
where $\delta >0$. This implies the decomposition (with the obvious notation) of the left-hand side of (5.16) into the sum $v_1(p,x)+v_2(p,x)$. The function $v_1(p,x)$ may be differentiated under the integral, if $x$ is in a small neighbourhood of $x^0$ (then $|x-\xi |$ remains separated from zero). After that, we set $x=x^0$, so that
\begin{equation}
\frac{\partial v_1(p,x^0)}{\partial x_j}=\int\limits_{\Omega_1} \frac{\partial^2 }{\partial x_j^2}\widetilde{Z}(p,x^0-\xi )\left[ [w_0(\xi )-w_0(x^0)\right] \,d\xi .
\end{equation}

Let $0<d<\delta /2$, $\widetilde{d}=(0,\ldots ,d,\ldots ,0)$ ($d$ is in the $j$-th place). Then
\begin{multline*}
\frac{1}d\left[ v_2(p,x^0+\widetilde{d})-v_2(p,x^0)\right] -\int\limits_{\Omega_2}\frac{\partial^2
\widetilde{Z}(p,x^0-\xi )}{\partial x_j^2}[w_0(\xi )-w_0(x^0)]\,d\xi
\\
=\frac{1}d\int\limits_{|x^0-\xi |\le 2d} \frac{\partial
\widetilde{Z}(p,x^0+\widetilde{d}-\xi )}{\partial x_j}[w_0(\xi )-w_0 (x^0)]\,d\xi
-\frac{1}d\int\limits_{|x^0-\xi |\le 2d} \frac{\partial
\widetilde{Z}(p,x^0-\xi )}{\partial x_j}[w_0(\xi )-w_0(x^0)]\,d\xi \\
-\int\limits_{|x^0-\xi |\le 2d} \frac{\partial^2
\widetilde{Z}(p,x^0-\xi )}{\partial x_j^2}[w_0(\xi )-w_0(x^0)]\,d\xi
\\
+\int\limits_{2d\le |x^0-\xi |\le \delta} \left\{ \frac{1}d\left[ \frac{\partial
\widetilde{Z}(p,x^0+\tilde{d}-\xi )}{\partial x_j}-\frac{\partial
\widetilde{Z}(p,x^0-\xi )}{\partial x_j}\right]-\frac{\partial^2
\widetilde{Z}(p,x^0-\xi )}{\partial x_j^2}\right\}[w_0(\xi )-w_0(x^0)]\,d\xi \\
\overset{\text{def}}{=}X_1+X_2+X_3+X_4.
\end{multline*}

We have
$$
|X_1|\le \frac{C}d\int\limits_{|x^0-\xi |\le 2d}\left| x^0+\widetilde{d}-\xi\right|^{-n+1+\gamma }\,d\xi
=\frac{C}d\int\limits_{|y|\le 2d}\left| y+\widetilde{d}\right|^{-n+1+\gamma }\,dy,
$$
and after the change $y=d\cdot \eta$ we find that $|X_1|\le Cd^\gamma \to 0$, as $d\to 0$. Similarly, $X_2\to 0$, as $d\to 0$. For $X_3$, we get from (5.12) that
$$
|X_3|\le C\int\limits_{|x^0-\xi |\le 2d}\left| x^0-\xi\right|^{-n+\gamma }\,d\xi =C_1d^\gamma \to 0.
$$

In $X_4$, by the Taylor formula, the expression in braces equals
$$
\frac{d}2\frac{\partial^3}{\partial x_j^3}\widetilde{Z}(p,x'-\xi ),\quad x'=x^0+\theta \widetilde{d},\ 0\le \theta \le 1,
$$
where
$$
|x'-\xi |\ge |\xi -x^0|-|x'-x^0|\ge |\xi -x^0|-d\ge \frac12 |\xi -x^0|,
$$
so that
$$
X_4\le Cd\int\limits_{2d\le |x^0-\xi |\le \delta}\left| x^0-\xi\right|^{-n+\gamma -1}\,d\xi =Cd^\gamma \int\limits_{2<|z|\le \delta d^{-1}}|z|^{-n-1+\gamma }\,dz \to 0,
$$
as $d\to 0$.

Hence,
$$
\frac{\partial v_2(p,x^0)}{\partial x_j}=\int\limits_{\Omega_2} \frac{\partial^2 }{\partial x_j^2}\widetilde{Z}(p,x^0-\xi )\left[ [w_0(\xi )-w_0(x^0)\right] \,d\xi .
$$
Taking into account (5.17) we find that, for any $p>0$, $x\in \Rn$,
\begin{equation}
\Delta_x\wpx =\int\limits_{\Rn}\Delta_x\left[ \widetilde{Z}(p,x-\xi )\right] [w_0(\xi )-w_0(x)]\,d\xi .
\end{equation}

For each $p>0$, the function $\Zpx$ is a fundamental solution of the equation $-\Delta u+p\K (p)u=0$ (it coincides, up to an easy change of variables, with the well-known fundamental solution for the equation $-\Delta u+u=0$; see Chapter 8 in \cite{Nik}). Therefore $\Delta \widetilde{Z}(p,x-\xi )=p\K (p)\widetilde{Z}(p,x-\xi )$, $x\ne \xi$. Substituting this in (5.18) we find that
$$
\Delta \wpx =p\K (p)\wpx -p\K (p)w_0(x)\int\limits_{\Rn}\widetilde{Z}(p,\xi )\,d\xi ,
$$
and it follows from (5.15) that $\wpx$ satisfies (5.1). $\qquad \blacksquare$

Note that a probabilistic representation of a fundamental solution of the problem (1.5) (and more general problems, with L\'evy generators instead of the Laplacian) was found in \cite{MS8,MNV}. The fundamental solutions were understood there as those of the equations obtained by applying the Laplace transform in time and the Fourier transform in spatial variables. For our situation, we obtained solutions, strong with respect to the variable $x$. To obtain classical solutions, one needs stronger assumptions regarding the function $\K (p)$. For example, it would be sufficient to assume its asymptotic properties found in \cite{K1} for the distributed order case (1.3). Then the investigation of strong solutions would repeat the reasoning from \cite{K1}.

Our uniqueness result for the problem (1.5) holds under much more general assumptions and is an immediate consequence of a deep result by E. E. Shnol (see Theorem 2.9 in \cite{CFKS}). Note that the notion of a LT-solution makes sense also for polynomially bounded solutions, that is such solutions $w(t,x)$ that $|w(t,x)|\le P(|x|)$ where $P$ is some polynomial independent of $t$. Instead of $(*)$, we make the following weaker assumption:

\begin{description}
\item[(**)] The function $k$ is nonnegative, locally integrable, nonzero on a set of positive measure, and its Laplace transform $\K (p)$ exists for all $p>0$.
\end{description}

\medskip
\begin{teo}
Let $(**)$ hold, and suppose that $w(t,x)$ is a polynomially bounded LT-solution of the problem (1.5) with $w_0(x)\equiv 0$. Then $w(t,x)\equiv 0$.
\end{teo}

\medskip
{\it Proof}. The Laplace transform $\wpx$ satisfies, for each $p>0$, the equation $\Delta \wpx =p\K (p)\wpx$. Thus $\wpx$ is a generalized eigenfunction of the operator $-\Delta$ on $L_2(\Rn )$ with the eigenvalue $-p\K (p)<0$. By Shnol's theorem, a nonzero polynomially bounded generalized eigenfunction is possible only if the eigenvalue belongs to the spectrum of $-\Delta$ equal to $[0,\infty )$. Therefore $\wpx \equiv 0$, so that $w(t,x)\equiv 0$. $\qquad \blacksquare$

\medskip
Theorem 4 can be extended to some equations with coefficients depending on $x$, for which Shnol's theorem can be applied; see \cite{CFKS,HS}.

\end{document}